\documentclass[submission]{FPSAC2019}

\usepackage{stmaryrd}

\allowdisplaybreaks

\newtheorem{thm}{Theorem}

\numberwithin{equation}{section}
\numberwithin{thm}{section}

\newcommand{\Z}{\mathbb{Z}}

\newcommand{\C}{\mathbb{C}}
\newcommand{\im }{\mathrm{i}}  
\newcommand{\dx }{\mathrm{d}}

\DeclareMathOperator*{\pf}{\mathrm{pf}}

\title[Edge asymptotics of skew Young diagrams]{New edge asymptotics of skew Young diagrams via free boundaries}

\author[D. Betea, J. Bouttier, P. Nejjar, M. Vuleti\'c]{Dan Betea\thanks{D.B. was partially supported by the German Research Foundation as part of the CRC 1060--B04 project.}\addressmark{1}, J\'er\'emie Bouttier\thanks{J.B. was partially supported by the ``Combinatoire \`a Paris'' project funded by the City of Paris, and by the grant ANR-14-CE25-0014 ``GRAAL''.}\addressmark{2, 3}, Peter Nejjar\thanks{P.N. was partially supported by ERC Advanced  Grant No. 338804 and ERC Starting Grant No. 716117.}\addressmark{4}, \and
Mirjana Vuleti\'c\addressmark{5}}

\address{
\addressmark{1}Institute for Applied Mathematics, University of Bonn, D-53115 Bonn, Germany \\ 
\addressmark{2}Institut de Physique Th\'eorique, Universit\'e Paris--Saclay, CEA, CNRS, F--91191 Gif--sur--Yvette, France \\
\addressmark{3}Univ Lyon, Ens de Lyon, Univ Claude Bernard, CNRS, Laboratoire de Physique, F-69342 Lyon, France \\
\addressmark{4}IST Austria, 3400 Klosterneuburg, Austria \\
\addressmark{5}Department of Mathematics, University of Massachusetts Boston, Boston, MA 02125, USA
}

\received{November 14, 2019}

\abstract{We study edge asymptotics of poissonized Plancherel-type measures on skew Young diagrams (integer partitions). These measures can be seen as generalizations of those studied by Baik--Deift--Johansson and Baik--Rains in resolving Ulam's problem on longest increasing subsequences of random permutations and the last passage percolation (corner growth) discrete versions thereof. Moreover they interpolate between said measures and the uniform measure on partitions. In the new KPZ-like $1/3$ exponent edge scaling limit with logarithmic corrections, we find new probability distributions generalizing the classical Tracy--Widom GUE, GOE and GSE distributions from the theory of random matrices.}

\usepackage[backend=bibtex]{biblatex}
\begin{document}

\maketitle

\section{Introduction and main results}

\textbf{Background.\ } The poissonized Plancherel \cite{bdj} and the discrete corner growth \cite{joh2} measures are two probability measures on integer partitions coming from the study of longest increasing subsequences of random permutations, respectively directed last passage percolation in an $N \times N$ square---LPP---models with iid geometric weights. Both (cf. op. cit.) are known to exhibit KPZ $N^{1/3}$ fluctuation behavior at the edge, with the Tracy--Widom GUE distribution~\cite{tw} from random matrix theory as the limiting distribution. Baik and Rains~\cite{br1, br2} have considered symmetrized versions of both, with similar results except now the limiting distributions are the Tracy--Widom GOE and GSE distributions~\cite{TW05} from random matrix theory---and some interpolating ones. In all cases the weight of a partition is proportional to the (possibly square of) number of (semi-) standard Young tableaux of said shape. Such measures are part of a large class of determinantal/pfaffian measures called Schur measures~\cite{oko}---see~\cite{br} for the pfaffian case.

In this note we generalize the above poissonized Plancherel and discrete Schur measures to skew Young diagrams---i.e.\,counting (semi-) standard Young tableaux of skew shape---coming from pairs (or triples) of partitions with free ends. In the $N^{1/3}$ scaling limit we find that the behavior of the first part of one (any) of the partitions is governed by new generalizations of the classical Tracy--Widom distributions. We finally give a directed LPP interpretation of our models in terms of directed polymers in a reflecting strip. We call these models \textit{LPP on a tie}.

\textbf{Outline.\ } Below we describe the main result; its relation to LPP in Section~\ref{sec:corner_growth}; the new distribution functions we obtain as limits in Section~\ref{sec:dist}; and a sketch of proof in Section~\ref{sec:proof}. Full details will appear elsewhere~\cite{BBNV2}.



\textbf{Main results.\ } For two integer partitions $\mu \subset \lambda$ (i.e.\,$\mu_i \leq \lambda_i\,\forall i$), let us denote by $f^{\lambda/\mu}$, respectively $\tilde{f}^{n,\,\lambda/\mu}$, the number of standard Young tableaux of skew shape $\lambda/\mu$ filled with numbers $1,2,\dots, |\lambda/\mu|:=|\lambda|-|\mu|$, respectively the number of semi-standard Young tableaux of shape $\lambda/\mu$ filled with $1, \dots, n$. Let $or(\lambda)$ (respectively $oc(\lambda)$) denote the number of \emph{odd rows} (respectively \emph{odd columns}) of a partition $\lambda$. Fix positive real parameters $u , q < 1$, $a_1, a_2, b_1, b_2, v \leq 1$\footnote{The restrictions on the $a_i,b_i$ can be somewhat relaxed.}, $\epsilon$\footnote{Here $u, v, a_i, b_j$ are boundary parameters, $\epsilon$ is a poissonization parameter keeping track of the size of the skew Young diagrams, and $q$ a geometrization parameter for the same purpose.}, and consider the following measures on pairs/triples of partitions:
\begin{equation} \label{eq:meas_def}
\begin{split}
  \mathbb{M}^{\nearrow,\,x} (\mu, \lambda) & \propto \Delta^x(\mu, \lambda) \cdot u^{|\mu|} \cdot \frac{\epsilon^{|\lambda/\mu|} f^{\lambda/\mu} } {|\lambda/\mu|!}, \\
  \mathbb{M}^{\nearrow \searrow,\,x} (\mu, \lambda, \nu) & \propto \Delta^x(\mu, \nu) \cdot u^{|\mu|} v^{|\nu|} \cdot \frac{\epsilon^{|\lambda/\mu| + |\lambda/\nu|} f^{\lambda/\mu} f^{\lambda/\nu} } {|\lambda/\mu|! \cdot |\lambda/\nu|!}, \\
  \widetilde{\mathbb{M}}^{\nearrow,\,x} (\mu, \lambda) & \propto \Delta^x(\mu, \lambda) \cdot u^{|\mu|} \cdot q^{|\lambda/\mu|} \tilde{f}^{n,\,\lambda/\mu}, \\
  \widetilde{\mathbb{M}}^{\nearrow \searrow,\,x} (\mu, \lambda, \nu) & \propto \Delta^x(\mu, \nu) \cdot u^{|\mu|} v^{|\nu|} \cdot q^{|\lambda/\mu| + |\lambda/\nu|} \tilde{f}^{n,\,\lambda/\mu} \tilde{f}^{n,\,\lambda/\nu}
\end{split}
\end{equation}
where $x \in \{aa, ab, bb, -\}$ is a \emph{boundary label} for the boundary partitions ($-$ stands for absence of boundary parameters); where 
\begin{equation} \label{eq:boundary_weights}
 \Delta^x(\mu, \lambda) = \begin{cases}
         a_1^{oc(\mu)} a_2^{oc(\lambda)}, & x = aa, \\
         a_1^{oc(\mu)} b_2^{or(\lambda)}, & x = ab, \\
         b_1^{or(\mu)} b_2^{or(\lambda)}, & x = bb, \\
         1, & x = -;
        \end{cases}
\end{equation}
and where the partition functions making each into a probability measure can be explicitly computed using the methods of~\cite{bbnv_17}.

They will be referred to as the \emph{upwards} (for $\nearrow$) and \emph{up-down} (for $\nearrow \searrow$) \textit{free-boundary poissonized Plancherel} (for $\mathbb{M}$) and \textit{free-boundary geometric corner growth} (for $\widetilde{\mathbb{M}}$) measures. 

The upwards measures $\mathbb{M}$ become the ``symmetrized'' poissonized Plancherel measures studied in~\cite{br1, br2} when $u = 0$. The up-down measures $\mathbb{M}$ become the classical poissonized Plancherel measure from~\cite{bdj} when $u = v = 0$. Moreover they all become the uniform measure on partitions when $\epsilon \to 0$. Thus they interpolate poissonized Plancherel $\leftrightarrow$ uniform\footnote{Another such interpolating measure has been studied in~\cite[Example 3.4]{B2007cyl} (for the bulk) and in~\cite[Theorem 1.1]{bb_18} (for the edge). It is an instance of the periodic Schur process, while the measures $\mathbb{M}$ considered in this paper are instances of the free boundary Schur process.}. A similar remark holds for the $\widetilde{\mathbb{M}}$ measures, replacing ``poissonized Plancherel'' with Johansson's corner growth~\cite{joh2} (for $\widetilde{\mathbb{M}}^{\nearrow \searrow,\,\cdots}$ when $u = v = 0$) and the respective symmetrized versions~\cite{br1} (for $\widetilde{\mathbb{M}}^{\nearrow,\,\cdots}$ when $u = 0$).

Our main result is a generalization of those in~\cite{bdj, br2}, and can be seen as a pfaffian analogue of~\cite[Theorem 1.1]{bb_18}. We concentrate on the large scale behavior of $\lambda_1$\footnote{We obtain essentially the same result if we consider $\mu$ or $\nu$ instead of $\lambda$.} as $\epsilon, n \to \infty$ while all other parameters go to 1 in a suitable critical regime.

\begin{thm} \label{thm:main_pl} 
Fix $\eta, \alpha_i, \beta_i, i=1,2$ positive reals. Let $M := \frac{\epsilon}{1-u^2} \to \infty$ and set $u = v = \exp(-\eta M^{-1/3})$ and $(a_i, b_i) = \left( u^{\alpha_i/\eta}, u^{\beta_i/\eta} \right)\ i=1,2$ all going to 1 as $M \to \infty$. (In particular, $\epsilon \sim M^{2/3} \to \infty$.) We have:
 \begin{equation}
  \begin{split}
  \lim_{M \to \infty} \mathbb{M}^{\nearrow,\,x} \left( \frac{\lambda_1 - 2M}{M^{1/3}} \leq s + \frac{1}{\eta} \log \frac{M^{1/3}}{\eta} \right) & = F^{1; d_x}(s), \\
  \lim_{M \to \infty} \mathbb{M}^{\nearrow \searrow,\,x} \left( \frac{\lambda_1 - 2M}{M^{1/3}} \leq s + \frac{1}{2 \eta} \log \frac{M^{1/3}}{2 \eta} \right) & = F^{2; d_x} (s)
  \end{split}
 \end{equation}
 where $x \in \{aa, ab, bb, -\}$, the distributions $F^{\,\cdots}$ are defined in Section~\ref{sec:dist}, and 
 \begin{equation} \label{eq:labels}
  d_x = \begin{cases}
         \alpha_1, \alpha_2; \eta, & x = aa, \\
         \alpha_1, -; \eta, & x = ab, \\
         \eta, & x = bb \text{ or } -.
        \end{cases}
 \end{equation}
\end{thm}

\begin{thm} \label{thm:main_disc}
  Fix $\eta, \alpha_i, \beta_i, i=1,2$ positive reals. As $n \to \infty$ ($n$ a positive integer), let $u = v = \exp(-\eta n^{-1/3})$, $(a_i, b_i) = \left( u^{\alpha_i/\eta}, u^{\beta_i/\eta} \right)\ i=1,2$ all going to 1 and set $q = 1-u^2 \to 0$. We have:
 \begin{equation}
  \begin{split}
  \lim_{n \to \infty} \widetilde{\mathbb{M}}^{\nearrow,\,x} \left( \frac{\lambda_1 - \chi n}{n^{1/3}} \leq s + \frac{1}{\eta} \log \frac{n^{1/3}}{\eta} \right) & = F^{1; d_x}(s), \\
  \lim_{n \to \infty} \widetilde{\mathbb{M}}^{\nearrow \searrow,\,x} \left( \frac{\lambda_1 - \chi n}{n^{1/3}} \leq s + \frac{1}{2 \eta} \log \frac{n^{1/3}}{2 \eta} \right) & = F^{2; d_x} (s)
  \end{split}
 \end{equation}
 where $\chi = 2 q \sum_{\ell \geq 0} \frac{u^{2 \ell}}{1-u^{2 \ell} q} \xrightarrow{n \to \infty} 2$, $x \in \{aa, ab, bb, -\}$, and $d_x$ is as in~\eqref{eq:labels}.
\end{thm}

Notice the logarithmic corrections in all cases. Also notice the unusual scaling $q = O(n^{-1/3})$ in Theorem~\ref{thm:main_disc}, different from the usual $q = O(n^{-1})$ \cite{joh4}.

We can also show convergence of the first $k$ parts of $\lambda$ to the first $k$ parts of the ensembles given by the corresponding kernels of Section~\ref{sec:dist}, a result in the spirit and generalizing those of~\cite{boo, joh4}. We omit the statement for brevity. 

We further emphasize we concentrate here on the new interesting ``crossover'' regime $u,v \to 1$. I.e.\,the case $u \to u_0 \in [0,1)$ leads, up to deterministic shift, to the same asymptotics as $u=0$---e.g.\ for the $x = -$ label, the limiting distributions are the Tracy--Widom $F_{\mathrm{GOE}}$ distribution in the case of the upwards measure and the Tracy--Widom $F_{\mathrm{GUE}}$ distribution for the up-down measure. We also omit this for brevity.

Finally, the new limiting distributions---defined in Section~\ref{sec:dist}---contain all the classical Tracy--Widom distributions as limits. 

\begin{thm}
  We have: $\lim_{\eta \to \infty} F^{1; \alpha_1, \alpha_2; \eta}(s) = F^{\boxslash}(s; \alpha_2),\ \lim_{\eta \to \infty} F^{2; \alpha_1, \alpha_2; \eta}(s) = F_{\mathrm{GUE}}(s)$ where $F_{\mathrm{GUE}}$ is the Tracy--Widom GUE distribution~\cite{tw} and $F^{\boxslash}(s; \alpha_2)$ is the Baik--Rains~\cite{br2} Tracy--Widom GOE/GSE~\cite{TW05} crossover---$F^{\boxslash}(s; 0) = F_{\mathrm{GOE}}(s)$ while $F^{\boxslash}(s; \infty) = F_{\mathrm{GSE}}(s)$. 
  \end{thm}

\section{A corner growth interpretation} \label{sec:corner_growth}

In this section we describe directed last passage percolation on a tie. More precisely, we show how the measures $\widetilde{\mathbb{M}}^{\nearrow \searrow,\,-}$ and $\widetilde{\mathbb{M}}^{\nearrow,\,-}$, and in particular the observables $\lambda_1$, come from certain LPP models on an infinite reflecting strip---the above mentioned \emph{tie models}. In the case of $u=v=0$ (for the up-down measure) and $u=0$ (for the upwards measure), the models become the usual LPP models of Johansson~\cite{joh2} and the Baik--Rains symmetrized versions~\cite{br1} respectively. For brevity, we will restrict the discussion to the measure $\widetilde{\mathbb{M}}^{\nearrow \searrow,\,-}$ and make remarks about the other one. 

\begin{figure}[!ht]
\begin{center}
 \includegraphics[scale=0.6]{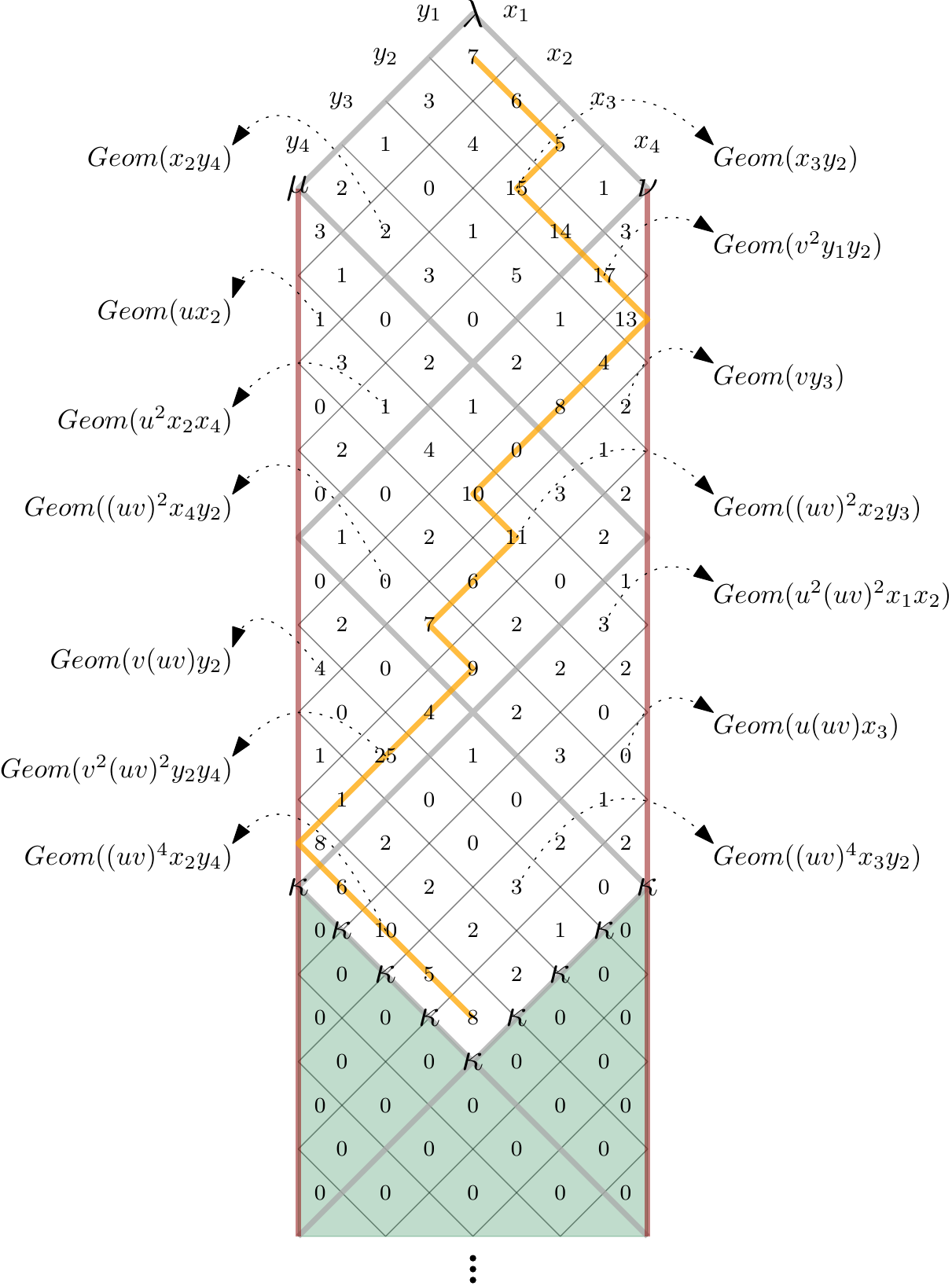} \qquad \includegraphics[scale=0.6]{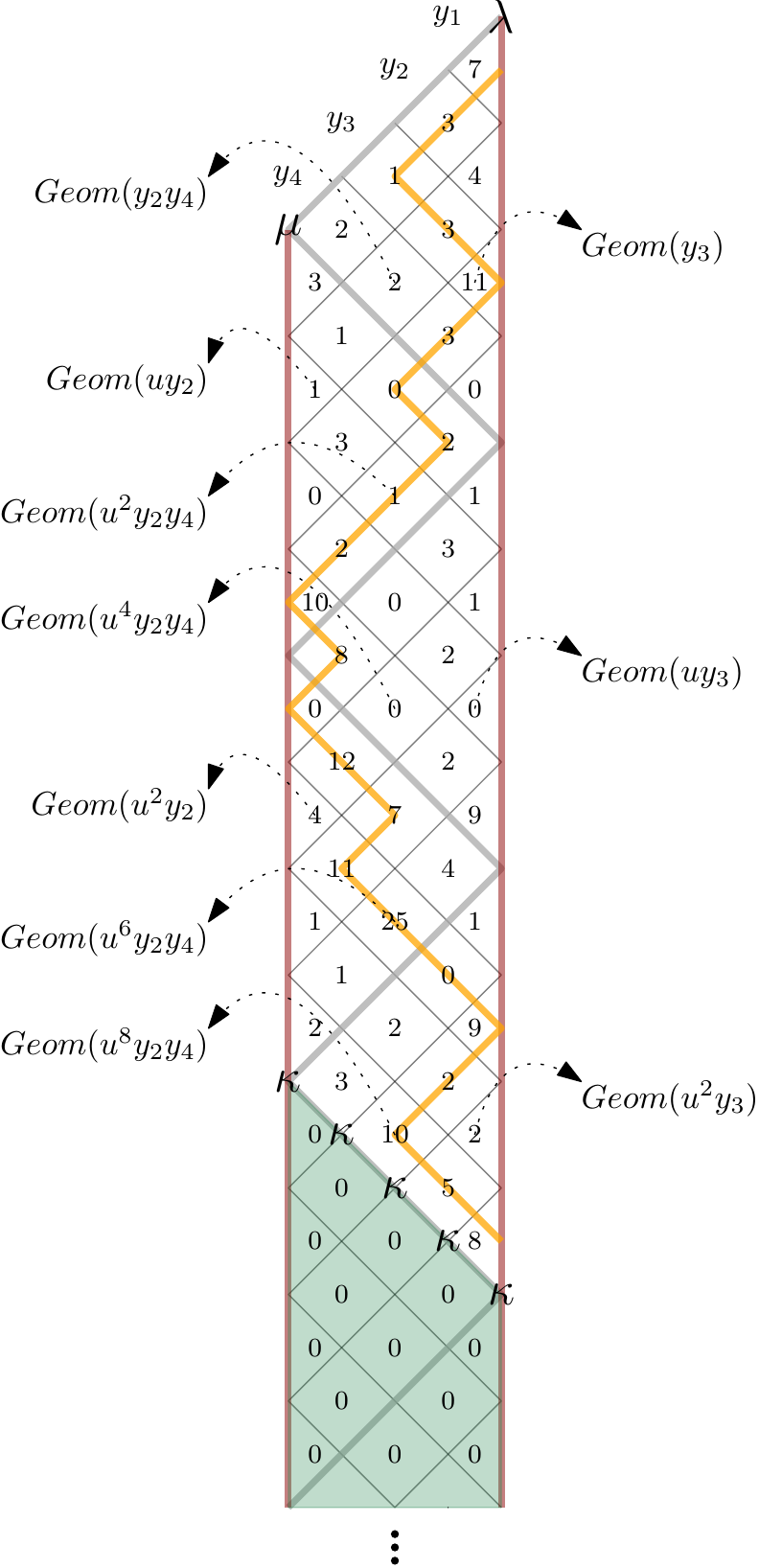}
\end{center}
 \caption{\footnotesize{The construction described in Section~\ref{sec:corner_growth}, with longest polymer in orange. $n=4$, $L = 199$ (left) and $L=130$ (right).}}
 \label{fig:fb_lpp}
\end{figure}

First we fix parameters $x_1, \dots, x_n, y_1, \dots, y_n$. Note in the end we take $x_i = y_i = q\ \forall i$. Suppose we have an infinite strip---or \emph{tie}---on the discrete square lattice constructed from big $n \times n$ adjacent squares and triangles like in Figure~\ref{fig:fb_lpp} (left). The strip has reflecting boundaries (red lines in fig. cit.). Each big square, sitting centrally in the strip, contains $n^2$ unit squares, and each big triangle $n(n-1)/2$ unit squares and $n$ unit triangles. In each unit square/triangle there is a geometric random variable $Geom(z)$\footnote{We say $X \sim Geom(z)$ if $\mathrm{Prob}(X=k) = (1-z) z^k, \ \forall k \geq 0$.}---independent from the others---of a certain parameter $z$ chosen as follows. Associate our $x_i, y_j$ parameters with the north-east (NE) and north-west (NW) boundaries of the strip (ends of the tie) as depicted in Figure~\ref{fig:fb_lpp} (left).  Pick a unit square from a big $n\times n$ square. The number inside has distribution $Geom((uv)^{2s} x_i y_j)$ where $s=0,1,2,\dots$ is the vertical position of the big square in the strip (starting from the top $s=0$) and to figure out $i,j$ send two rays of light from said square to the top, one in the NE and the other in the NW directions. The rays reflect off each boundary. They will intersect the top NE and NW borders in an $x_i$ and $y_j$ parameter respectively, which are our sought variables. For a unit square inside a big $n \times n$ triangle, the number inside is either $Geom(u^2 (uv)^{2s} x_i x_j)$ or $Geom(v^2 (uv)^{2s} y_i y_j)$ with $s$ and determining $i,j$ as above. Note in this case, the two rays will hit the same top boundary, either labeled $x$ or $y$. Finally, in the unit boundary triangles, the numbers inside are $Geom(u (uv)^s x_i)$ or $Geom(v (uv)^s y_i)$) with $s$ and determining $i$ as above. See Figure~\ref{fig:fb_lpp} (left) for precise examples. 

By Borel--Cantelli, almost surely only finitely many numbers in this strip will be non-zero---say those outside the green area in Figure~\ref{fig:fb_lpp} (left). Look at the longest polymer (path) with south-east (SE) or south-west (SW) steps starting from the top unit square in the strip and going down, reflected by the two boundaries if need be. Here by \emph{length} we mean the sum of the integers encountered by the path. Call this length $L$. It equals $199$ in Figure~\ref{fig:fb_lpp} (left). 

Pick now a uniformly distributed partition $\kappa$ with parameter $uv$---i.e.\,$\mathrm{Prob}(\kappa) = (uv; uv)_{\infty} (uv)^{|\kappa|}$---and place it at the \emph{nodes} of the south-eastern and south-western-most bottom boundaries separating the infinite region of 0's inside the strip---without loss of generality it can be positioned on the SE and SW sides of a big square. Using the Fomin growth rule description of the Robinson--Schensted--Knuth (RSK) correspondence---like described in e.g.\,\cite{bbbccv}---inductively ``flip'', starting from the bottom, every unit square and triangle to produce, from three (or two for the triangles) partitions on its boundary and the integer inside, a fourth partition---placed at the top vertex. Call the final partition sitting at the top of the tie $\lambda$, and $\mu, \nu$ the ones sitting at the top ends of the reflecting boundaries. The properties of the RSK correspondence imply $\mathrm{Prob}(\mu, \lambda, \nu) \propto \widetilde{\mathbb{M}}^{\nearrow \searrow,\,-} (\lambda, \mu, \nu)$ (upon taking $x_i = y_i = q \ \forall i$). Greene's theorem~\cite{gre} yields $L+\kappa_1=\lambda_1$. We thus obtain:

\begin{thm} \label{thm:corner_growth}
 It holds that $\mathrm{Prob}(L + \kappa_1 \leq k) = \widetilde{\mathbb{M}}^{\nearrow \searrow,\,-} (\lambda_1 \leq k)$.
\end{thm}

The corner growth construction corresponding to the measure $\widetilde{\mathbb{M}}^{\nearrow,\,-}$ is depicted in Figure~\ref{fig:fb_lpp} (right; $L=130$; there are no big squares and hence no $x$ parameters; $v=1$). Theorem~\ref{thm:corner_growth} holds mutatis-mutandis. A similar construction holds for the $\mathbb{M}^{\nearrow,\,-}, \mathbb{M}^{\nearrow \searrow,\,-}$ measures if one replaces the geometric random variables with appropriate Poisson point processes. For the measures with boundary $a$ and/or $b$ parameters, a more involved construction generalizing that of~\cite{br1} exists. E.g.\,for the $\widetilde{\mathbb{M}}^{\nearrow \searrow,\,aa}$ measure, one just puts another distribution in the unit triangles of Figure~\ref{fig:fb_lpp} (left): $Geom(a_1 u (uv)^s x_i)$ respectively $Geom(a_2 v (uv)^s y_i)$ for the two boundaries. For a $b$ parameter---say on the right boundary---one replaces $Geom(a_2 v (uv)^s y_i)$ with $Geom^{(b_2)} (v (uv)^s y_i)$ where $X \sim Geom^{(b)} (z)$ if $\mathrm{Prob} (X = k) = (1-z^2)(1+bz)^{-1} b^{k \ \mathrm{mod}\ 2} z^{k}\ \forall k\geq 1$. Because of the mod 2, one heuristically sees that the $b$ parameters do not matter in the limits of Theorem~\ref{thm:main_disc}, whereas the $a$'s do and can be taken as ``strengths'' of the respective boundaries.

\section{Definition of distribution functions} \label{sec:dist}

Fix parameters $\alpha_1, \alpha_2, \eta > 0$. Our new limit distributions are defined as Fredholm pfaffians. For $k=1,2$ let 
\begin{equation} 
\begin{split}
F^{k; \alpha_1, \alpha_2; \eta} (s) &:= \pf \left( J - A^{k; \alpha_1, \alpha_2; \eta} \right)_{L^2\left( s + \frac{\log2}{k \cdot \eta}, \infty \right)},\\
F^{k; \alpha_1, -; \eta} &:= \lim_{\alpha_2 \to 0} F^{k; \alpha_1, \alpha_2; \eta}, \qquad \qquad \qquad \qquad F^{k; \eta} := \lim_{\alpha_1 \to 0} F^{k; \alpha_1, -; \eta}
\end{split}
\end{equation}
where $J$ is the anti-symmetric kernel $J(x,y) = \delta_{x, y} \left( \begin{array}{rr} 0 & 1 \\ -1 & 0 \end{array} \right)$ and where the $A$ kernels (operators) are defined as follows. Let $\tau, \tau'$ satisfy $0 < \tau, \tau' < \frac{1}{6}\min (\eta, \alpha_1, \alpha_2)$. In case $\alpha_1$ (or $\alpha_2$ or both) is zero, the corresponding factor is absent from the minimum conditions. Define the following pairs of contours, oriented bottom-to-top: $(C^{1,1}_{\zeta}, C^{1,1}_{\omega}) = (\tau + \im \mathbb{R}, \tau' + \im \mathbb{R}), (C^{1,2}_{\zeta}, C^{1,2}_{\omega}) = (\tau + \im \mathbb{R}, -\tau' + \im \mathbb{R}), (C^{2,2}_{\zeta}, C^{2,2}_{\omega}) = (-\tau + \im \mathbb{R}, -\tau' + \im \mathbb{R}).$ Define these auxiliary products of Euler Gamma functions\footnote{We use the notation $\Gamma(a, b, c, \dots) = \Gamma(a) \Gamma(b) \Gamma(c) \cdots$.}: 
\begin{equation}
 \begin{split}
  \gamma^{(1)} (\zeta) := \frac{ \Gamma\left(\frac{1}{2} + \frac{\alpha_1 - \zeta}{2 \eta}, 1 + \frac{\alpha_2 - \zeta}{2 \eta}\right) }{\Gamma\left(\frac{1}{2} + \frac{\alpha_1 + \zeta}{2 \eta}, \frac{\alpha_2 + \zeta}{2 \eta}\right)}, \qquad \gamma^{(2)}(\zeta) &:= \frac{ \Gamma\left(\frac{1}{4} + \frac{\alpha_1 - \zeta}{4 \eta}, \frac{3}{4} + \frac{\alpha_2 - \zeta}{4 \eta}\right) }{\Gamma\left(\frac{3}{4} + \frac{\alpha_1 + \zeta}{4 \eta}, \frac{1}{4} + \frac{\alpha_2 + \zeta}{4 \eta}\right)}.
 \end{split}
\end{equation}
The desired $2 \times 2$ \textit{hypergeometric Airy} matrix kernels are below. The $\zeta$ integral is always over $C^{i,j}_{\zeta}$, the $\omega$ always over $C^{i,j}_{\omega}$ for corresponding $i,j$, $\dx_{\zeta \omega} := \dx \zeta \dx \omega (2 \pi \im)^{-2}$, and $\dx_{\zeta} := \dx \zeta (2 \pi \im)^{-1}$.
\begin{align*}
A^{1; \alpha_1, \alpha_2; \eta}_{1,1}(x, y) =& \int \int  \Gamma \left( \frac{\zeta}{\eta}, \frac{\omega}{\eta} \right) \gamma^{(1)} (\zeta) \gamma^{(1)} (\omega)  \frac{\sin \frac{\pi(\zeta-\omega)}{2 \eta} } {\sin \frac{\pi(\zeta+\omega)}{2 \eta}}    e^{ \frac{\zeta^3}{3} - x \zeta + \frac{\omega^3}{3} - y \omega } \dx_{\zeta \omega}, \\
A^{1; \alpha_1, \alpha_2; \eta}_{1,2}(x, y) =& \int \int \Gamma \left( \frac{\zeta}{\eta}, 1 - \frac{\omega}{\eta} \right) \frac{\gamma^{(1)} (\zeta)}{\gamma^{(1)} (\omega) } \frac{\sin \frac{\pi(\zeta+\omega)}{2 \eta} } {\sin \frac{\pi(\zeta-\omega)}{2 \eta}}    e^{ \frac{\zeta^3}{3} - x \zeta - \frac{\omega^3}{3} + y \omega }  \frac{\dx_{\zeta \omega}}{2 \eta}, \tag{\stepcounter{equation}\theequation} \\
A^{1; \alpha_1, \alpha_2; \eta}_{2,2}(x, y) =& \int \int \Gamma \left(1 - \frac{\zeta}{\eta}, 1 - \frac{\omega}{\eta} \right) \frac{1} {\gamma^{(1)} (\zeta) \gamma^{(1)} (\omega)}  \frac{\sin \frac{\pi(\zeta-\omega)}{2 \eta} } {\sin \frac{\pi(\zeta+\omega)}{2 \eta}}    e^{ - \frac{\zeta^3}{3} + x \zeta - \frac{\omega^3}{3} + y \omega } \frac{\dx_{\zeta \omega}}{4 \eta^2} \\
& - \mathrm{sgn}(x-y)
\end{align*}
and
\begin{align*}
A^{2; \alpha_1, \alpha_2; \eta}_{1,1}(x, y) &= \int \int \Gamma \left(\frac{1}{2} + \frac{\zeta}{2 \eta}, \frac{1}{2} + \frac{\omega}{2 \eta} \right) \gamma^{(2)} (\zeta) \gamma^{(2)} (\omega) \frac{\sin \frac{\pi(\zeta-\omega)}{4 \eta} } {\cos \frac{\pi(\zeta+\omega)}{4 \eta}}    e^{ \frac{\zeta^3}{3} - x \zeta + \frac{\omega^3}{3} - y \omega }  \frac{\dx_{\zeta \omega}}{4 \eta}, \\
A^{2; \alpha_1, \alpha_2; \eta}_{1,2}(x, y) &= \int \int \Gamma \left(\frac{1}{2} + \frac{\zeta}{2 \eta}, \frac{1}{2} - \frac{\omega}{2 \eta} \right) \frac{\gamma^{(2)} (\zeta)}{\gamma^{(2)} (\omega)} \frac{\cos \frac{\pi(\zeta+\omega)}{4 \eta} } {\sin \frac{\pi(\zeta-\omega)}{4 \eta}}   e^{ \frac{\zeta^3}{3} - x \zeta - \frac{\omega^3}{3} + y \omega }  \frac{\dx_{\zeta \omega}}{4 \eta}, \tag{\stepcounter{equation}\theequation} \\
A^{2; \alpha_1, \alpha_2; \eta}_{2,2}(x, y) &= \int \int  \Gamma \left(\frac{1}{2} - \frac{\zeta}{2 \eta}, \frac{1}{2} - \frac{\omega}{2 \eta} \right) \frac{1}{\gamma^{(2)} (\zeta) \gamma^{(2)} (\omega)} \frac{\sin \frac{\pi(\zeta-\omega)}{4 \eta} } {\cos \frac{\pi(\zeta+\omega)}{4 \eta}}    e^{ -\frac{\zeta^3}{3} + x \zeta - \frac{\omega^3}{3} + y \omega } \frac{\dx_{\zeta \omega}}{4 \eta}
\end{align*}
with the remark that everywhere $A_{2, 1}^{\,\cdots} (x, y) := - A_{1, 2}^{\,\cdots} (y, x)$. 

We note that, upon using the Gamma duplication formula, we can write, for $k=1,2$, 
\begin{equation}
 F^{k; \eta} \left( s-\frac{\log 2}{k \cdot \eta} \right) = \pf \left( J - A^{k; \eta} \right)_{L^2( s, \infty )}
\end{equation}
where the kernels $A^{k; \eta}$ have a much simpler form:
\begin{align*}
A^{1; \eta}_{1,1}(x, y) =& \int \int  \Gamma \left( 1 - \frac{\zeta}{\eta}, 1 - \frac{\omega}{\eta} \right) \frac{\sin \frac{\pi(\zeta-\omega)}{2 \eta} } {\sin \frac{\pi(\zeta+\omega)}{2 \eta}}   e^{ \frac{\zeta^3}{3} - x \zeta + \frac{\omega^3}{3} - y \omega }   \frac{\dx_{\zeta \omega}}{4}, \\
A^{1; \eta}_{1,2}(x, y) =& \int \int \Gamma \left( 1 - \frac{\zeta}{\eta}, \frac{\omega}{\eta} \right) \frac{\sin \frac{\pi(\zeta+\omega)}{2 \eta} } {\sin \frac{\pi(\zeta-\omega)}{2 \eta}}    e^{ \frac{\zeta^3}{3} - x \zeta - \frac{\omega^3}{3} + y \omega }   \frac{\dx_{\zeta \omega}}{2 \eta}, \tag{\stepcounter{equation}\theequation}\\
A^{1; \eta}_{2,2}(x, y) =& \int \int \Gamma \left( \frac{\zeta}{\eta}, \frac{\omega}{\eta} \right) \frac{\sin \frac{\pi(\zeta-\omega)}{2 \eta} } {\sin \frac{\pi(\zeta+\omega)}{2 \eta}}    e^{ - \frac{\zeta^3}{3} + x \zeta - \frac{\omega^3}{3} + y \omega }   \frac{\dx_{\zeta \omega}}{\eta^2} \\
& + \int \Gamma\left( \frac{\zeta}{\eta} \right) e^{- \frac{\zeta^3}{3} + x \zeta} \frac{\dx_\zeta}{\eta} - \int \Gamma\left( \frac{\omega}{\eta} \right) e^{- \frac{\omega^3}{3} + y \omega} \frac{\dx_\omega}{\eta} - \mathrm{sgn}(x-y)
\end{align*}
and
\begin{equation}
\begin{split}
A^{2; \eta}_{1,1}(x, y) =& \int \int  \Gamma \left( \frac{1}{2} - \frac{\zeta}{2 \eta}, \frac{1}{2} - \frac{\omega}{2 \eta} \right) \frac{\sin \frac{\pi(\zeta-\omega)}{4 \eta} } {\cos \frac{\pi(\zeta+\omega)}{4 \eta}}   e^{ \frac{\zeta^3}{3} - x \zeta + \frac{\omega^3}{3} - y \omega } \frac{\dx_{\zeta \omega}}{4 \eta}, \\
A^{2; \eta}_{1,2}(x, y) =& \int \int \Gamma \left( \frac{1}{2} - \frac{\zeta}{2 \eta}, \frac{1}{2} + \frac{\omega}{2 \eta} \right) \frac{\cos \frac{\pi(\zeta+\omega)}{4 \eta} } {\sin \frac{\pi(\zeta-\omega)}{4 \eta}}    e^{ \frac{\zeta^3}{3} - x \zeta - \frac{\omega^3}{3} + y \omega } \frac{\dx_{\zeta \omega}}{4 \eta}, \\
A^{2; \eta}_{2,2}(x, y) =& \int \int \Gamma \left( \frac{1}{2} + \frac{\zeta}{2 \eta}, \frac{1}{2} + \frac{\omega}{2 \eta} \right) \frac{\sin \frac{\pi(\zeta-\omega)}{4 \eta} } {\cos \frac{\pi(\zeta+\omega)}{4 \eta}}    e^{ - \frac{\zeta^3}{3} + x \zeta - \frac{\omega^3}{3} + y \omega } \frac{\dx_{\zeta \omega}}{4 \eta}
  \end{split}
\end{equation}
where again $A_{2, 1}^{\,\cdots} (x, y) := - A_{1, 2}^{\,\cdots} (y, x)$ and the contours are as before with one important exception: in the case of $A^{1;\,\eta}_{1,2}$, the $\omega$ contour $C^{1,2}_{\omega}$ passes locally \emph{to the right} of 0, but is otherwise as stated (to account for the pole at $\omega = 0$ whose residue was not taken). This happens because---up to inessential conjugation and taking the appropriate residues in $A^{1;\,\cdots}_{2,2}$ as $\alpha_1, \alpha_2 \to 0$---we have $A^{k; 0, 0; \eta} \left(x+\frac{2 \log 2}{k \cdot \eta}, y+\frac{2 \log 2}{k \cdot \eta} \right) = A^{k; \eta} (x,y),\ k=1,2$.

\section{Sketch of proof} \label{sec:proof}

The first tool we use to prove our results is the Schur measure with two free boundaries from~\cite{bbnv_17}. Recall the definition of skew Schur functions evaluated at a \textit{specialization} $\rho$ via the Jacobi--Trudi formula: $s_{\lambda/\mu}(\rho) = \det_{1 \leq i,j \leq n} h_{\lambda_i - i - (\mu_j - j)}(\rho)$ for $n$ large enough. Here a specialization $\rho$ is just a sequence of numbers $(h_n(\rho))_{n \geq 0}$---its values on the complete symmetric functions---assembled into the generating series $H(\rho; z) := \sum_{n \geq 0} h_{n}(\rho) z^n$. 

On sequences of partitions $\mu \subset \lambda \supset \nu$, consider the weights
\begin{equation} \label{eq:fbschur_weight}
  \mathcal{W}^{x}(\lambda, \mu, \nu) := \Delta^{x} (\mu, \nu) \cdot u^{|\mu|} v^{|\nu|} \cdot s_{\lambda / \mu}(\rho^+)
    s_{\lambda / \nu} ( \rho^- )
\end{equation}
with $x \in \{aa, ab, bb, -\}$, $\rho^{\pm}$ two specializations and $\Delta^x$ as in~\eqref{eq:boundary_weights}---note after proper normalization~\cite{bbnv_17}, these become probability measures. We observe our original measures~\eqref{eq:meas_def} are of the form given in equation~\eqref{eq:fbschur_weight} as follows. For the up-down $\nearrow \searrow$ measures, $\rho^+ = \rho^-$ and both are the poissonized Plancherel $pl_{\epsilon}$ specialization for $\mathbb{M}$ and $\rho_q$ (a specialization in $n$ variables all $=q$) for $\widetilde{\mathbb{M}}$ respectively; for the upwards measures $\nearrow$, set $v=1$ and $\rho^+ = pl_{\epsilon}$ for $\mathbb{M}$, $=\rho_q$ for $\widetilde{\mathbb{M}}$ while in both cases $\rho^- = 0$, the empty specialization; recall $\rho_{\epsilon}, \rho_q$ are defined by $H(pl_{\epsilon}; z) = \exp(\epsilon z), H(\rho_q; z) = (1-qz)^{-n}$; and finally note parameters $a_1, a_2, b_1, b_2, u, v, q$ satisfy the inequalities from the Introduction. The reason for the above is as follows: Schur functions, specialized in variables, are generating series for semi-standard Young tableaux. As such, $s_{\lambda / \mu}(\rho_q) = s_{\lambda / \mu}(q, \dots, q) = q^{|\lambda/\mu|} \tilde{f}^{n,\,\lambda/\mu}$ (there are $n$ $q$'s inside the parentheses) and a limiting argument shows $s_{\lambda / \mu} (pl_{\epsilon}) = \epsilon^{|\lambda / \,u|} \frac{f^{\lambda/\mu}}{|\lambda/\mu|!}$.

Define the \emph{point configuration} associated with $\lambda$ coming from the above measures as $\mathfrak{S}^x(\lambda) := \{ \lambda_j-j+1/2 \} \subset \Z+1/2$. It is a simple point process. Consider the \emph{shifted} version $\mathfrak{S}^x_t(\lambda):=\mathfrak{S}^x(\lambda) + 2D_t$ where for $t$ a parameter, $D_t$ is an integer independent of everything else having \emph{theta} distribution
$\mathrm{Prob}(D_t=d) = \frac{t^{2d} (uv)^{2d^2}}{\theta_3(t^2;(uv)^4)}$. A slightly more general version of~\cite[Theorem 2.5]{bbnv_17} shows that, for $x \in \{aa, ab, bb, -\}$, the shifted point process
$\mathfrak{S}^x_t(\lambda)$ is pfaffian with $2 \times 2$ matrix correlation kernel $K^x$ of the form:
  \begin{equation} \label{eq:Kint}
    \begin{split}
      K_{1,1}^x(k; k') &= \frac{1}{(2\im \pi)^2}
      \oint_{|z|=r} \frac{\dx z}{z^{k+1}} \oint_{|w|=r'} \frac{\dx w}{w^{k'+1}}
      F(z) F(w) \kappa^x_{1,1}(z,w), \\
      K_{1,2}^x(k; k') &= - K_{2,1}^x(k';k)
       = \frac{1}{(2\im \pi)^2}
      \oint_{|z|=r} \frac{\dx z}{z^{k+1}} \oint_{|w|=r'} \frac{\dx w}{w^{-k'+1}}
      \frac{F(z)}{F(w)} \kappa^x_{1,2}(z,w), \\
      K_{2,2}^x(k; k') &= \frac{1}{(2\im \pi)^2}
      \oint_{|z|=r} \frac{\dx z}{z^{-k+1}} \oint_{|w|=r'} \frac{\dx w}{w^{-k'+1}}
      \frac{1}{F(z) F(w)} \kappa^x_{2,2}(z,w)
    \end{split}
  \end{equation}
with $F(z) = \frac{ H(\rho^+;z)}{
        H(\rho^-;z^{-1})} 
     \prod_{n \geq 1} \frac{H(u^{2n} v^{2n-2} \rho^-;  z) H(u^{2n} v^{2n} \rho^+; z)}{
       H(u^{2n-2} v^{2n} \rho^+; z^{-1}) H(u^{2n} v^{2n} \rho^-; z^{-1})}$;
with $\kappa_{1,1}^x(z,w)$, $\kappa_{1,2}^x(z,w)$ and $\kappa_{2,2}^x(z,w)$ given\footnote{We use the notation $(x; q)_\infty := \prod_{\ell \geq 0} (1-x q^\ell)$, $\theta_q (x) := (x; q)_{\infty} (q/x; q)_{\infty}$, $\theta_3 (z; q) := (q;q)_{\infty} \theta_q(-z \sqrt{q})$  and $(a, b, \dots; q)_\infty := (a; q)_\infty \cdot (b; q)_\infty \cdots$.} respectively by
\begin{align*} 
&\frac{v^2 }{t z^{1/2} w^{3/2}} \cdot \frac{((uv)^2; (uv)^2)^2_{\infty}}  {(uz, uw, - \frac{v}{z}, - \frac{v}{w}; uv)_{\infty}} \cdot  \frac{\theta_{(uv)^2}(\frac{w}{z})}{\theta_{(uv)^2} (u^2 zw)} \cdot \frac{\theta_3 \left( ( \frac{t zw}{v^2})^2 ; (uv)^4 \right)}{\theta_3(t^2;(uv)^4)} \cdot g^{x}(z) g^x(w), \\
&\frac{w^{1/2}}{z^{1/2}} \cdot \frac{((uv)^2; (uv)^2)^2_{\infty}}{(uz, -uw,- \frac{v}{z}, \frac{v}{w}; uv)_{\infty}} \cdot \frac{\theta_{(uv)^2}(u^2 zw)}{\theta_{(uv)^2} (\frac{w}{z})} \cdot \frac{\theta_3 \left( ( \frac{t z}{w})^2 ; (uv)^4 \right)}{\theta_3(t^2;(uv)^4)} \cdot \frac{g^{x}(z)}{g^x(w)}, \tag{\stepcounter{equation}\theequation} \\
&\frac{t v^2}{z^{1/2} w^{3/2}} \cdot \frac{((uv)^2; (uv)^2)^2_{\infty}}{(-uz, -uw,\frac{v}{z}, \frac{v}{w}; uv)_{\infty}} \cdot \frac{\theta_{(uv)^2}(\frac{w}{z})}{\theta_{(uv)^2} (u^2 zw)} \cdot \frac{\theta_3 \left( ( \frac{t v^2}{zw})^2 ; (uv)^4 \right)}{\theta_3(t^2;(uv)^4)} \cdot \frac{1}{g^{x}(z) g^x(w)}
\end{align*}
where $g^x(z) = $
\begin{equation}
    \begin{cases}
        \frac{(uz; uv)_{\infty}}{\left( \frac{v}{z}; uv\right)_{\infty}} \cdot \frac{\left(\frac{a_1 u v^2}{z}, \frac{a_2 v}{z}; (uv)^2 \right)_{\infty}}{\left(a_1 u z, a_2 u^2 v z; (uv)^2 \right)_{\infty}} \cdot \frac{1}{(a_1 a_2 uv; (uv)^2)_{\infty} (-uv; uv)_{\infty}}, & x = aa, \\
        \frac{\left( -\frac{v}{z}; uv \right)_{\infty}} {( -uz; uv)_{\infty}} \cdot \frac{ (- b_1 uz, -b_2 u^2 vz; (uv)^2 )_{\infty} } {\left( - \frac{b_1 u v^2}{z}, - \frac{b_2 v}{z}; (uv)^2 \right)_{\infty}} \cdot \frac{1}{(b_1 b_2 uv; (uv)^2)_{\infty} (-uv; uv)_{\infty}}, & x = bb, \\
        \frac{\left( uz, -\frac{v}{z}, \frac{a_1 u v^2}{z}, -b_2 u^2 v z, -uv, -a_1 b_2 uv, a_1 (uv)^2, b_2 (uv)^2; (uv)^2\right)_{\infty}}{\left( -u^2 v z, \frac{u v^2}{z}, a_1 u z, - \frac{b_2 v}{z}, -a_1 uv, -b_2 uv; (uv)^2\right)_{\infty}} \cdot \frac{1}{(-uv, -uv, a_1 uv, b_2 uv; uv)_{\infty}}, & x = ab, \\
        1, & x = -;
    \end{cases}
\end{equation}
where $r, r'$ satisfy $\max(v, q) < r, r' < \min(q^{-1}, u^{-1})$ ($q$ does not appear for measures $\mathbb{M}$) and $r' < r$ for $K^{x}_{1,2}$; and where the choice of $\rho^\pm$ was described in the previous paragraph. To that point, we have $F(z) = \exp \left( \frac{\epsilon}{1-u^2} (z - z^{-1}) \right)$ (for all  $\mathbb{M}$) and $F(z) = \frac{(q/z; u^2)^n_{\infty}}{(q z; u^2)^n_{\infty}}$ (for all $\widetilde{\mathbb{M}}$). In this shifted process, the distribution of $\lambda_1$, in each case, is a Fredholm pfaffian for its respective kernel. 

Using steepest descent asymptotic analysis, we can show that the the latter kernels converge to the ones defining our distributions in the limits prescribed in Theorems~\ref{thm:main_pl} and~\ref{thm:main_disc}. We write $F(z) z^{-k}$ as $\exp(M S(z))$ (for $\mathbb{M}$) and $\exp(n S(z))$ (for $\widetilde{\mathbb{M}}$) in the scaling prescribed with $S$ having a double critical point at $z=1$; we thus scale the integration variables $(z, w) = (\exp(\zeta M^{-1/3}), \exp(\omega M^{-1/3}))$ (replace $M$ by $n$ for $\widetilde{\mathbb{M}}$) so that only an $M^{1/3}$ (or $n^{1/3}$) neighborhood of $1$ makes a non-zero contribution---see~\cite[Section 5]{bb_18} for a related computation. For $\kappa$, one uses the estimate
\begin{equation}
    \log (q^c; q)_{\infty} = -\frac{\pi^2}{6} r^{-1} + \left( \frac{1}{2} - c \right) \log r + \frac{1}{2} \log(2 \pi) - \log \Gamma(c) + O(r)
\end{equation} 
for $q = e^{-r}, r \to 0^+$ ($c \in \C-\Z_{\leq 0}$), along with $(-x; q)_{\infty} = \frac{(x^2; q^2)_{\infty}}{(x; q)_{\infty}}$. The logarithmic shifts in our results come from the $\log r$ terms above.

This together yields the result modulo the shift we have made. It ($2 D_t$) can be removed in the end\footnote{The argument of~\cite[Section 5]{bb_18} applies.} as $2 D_t / M^{1/3}$ (or $2 D_t / n^{1/3}$) converges to 0 in probability. Moreover, a slightly improved version of the above argument can be used to show convergence of multi-point correlation functions to the desired ensembles---one again has to remove the shift at the end. This finishes the proof.

\acknowledgements{D.B. is especially grateful to Patrik Ferrari for suggesting simplifications in Section~\ref{sec:dist} and to Alessandra Occelli for suggesting the name for the models of Section~\ref{sec:corner_growth}.}

\printbibliography
\end{document}